\documentclass{amsproc}
\usepackage{graphicx}
\usepackage{amscd}

\newtheorem{theorem}{Theorem}[section]
\newtheorem{corollary}[theorem]{Corollary}
\newtheorem{lemma}[theorem]{Lemma}
\newtheorem{proposition}[theorem]{Proposition}
\theoremstyle{definition}
\newtheorem{definition}[theorem]{Definition}

\theoremstyle{remark}

\numberwithin{equation}{section}

\begin{document}
\title[Fixed Points and Selections]{Fixed points and selections of set-valued maps on spaces with convexity.}
\author{Peter Saveliev}
\address{Department of Mathematics, University of Illinois, 1409 West Green Street,
Urbana, IL 61801}
\email{saveliev@math.uiuc.edu}
\date{December 26, 1997}
\subjclass{47H04, 47H10, 52A01, 54C65, 54H25}
\keywords{fixed point, continuous selection, generalized convexity}

\begin{abstract}
We provide theorems extending both Kakutani and Browder fixed points
theorems for multivalued maps on topological vector spaces, as well as some
selection theorems. For this purpose we introduce convex structures more
general than those of locally convex and non-locally convex topological
vector spaces or generalized convexity structures due to Michael, Van de Vel
and Horvath.
\end{abstract}

\maketitle

\section{Introduction.}

The origin of our fixed point theorem is the following two classical results
due to Kakutani-Ky Fan-Glicksberg \cite{Ka,Fa,Gl} and Browder \cite{Br}
respectively (see also \cite{Zeidler}).

\begin{theorem}[Kakutani Fixed-Point Theorem]
\label{Ka}Let $X$ be a nonempty convex compact subset of a locally convex
Hausdorff topological vector space, and let $F:X\rightarrow X$ be an u.s.c.
multifunction with nonempty convex closed images. Then $F$ has a fixed point.
\end{theorem}

\begin{theorem}[Browder Fixed-Point Theorem]
\label{Br}Let $X$ be a nonempty convex compact subset of a Hausdorff
topological vector space, and let $G:X\rightarrow X$ be a multifunction with
nonempty convex images and preimages relatively open in $X$. Then $G$ has a
fixed point.
\end{theorem}

Similarly, our selection theorem unites the following two results due to
Michael \cite{Mich2} and Browder \cite{Br} (see also \cite{Zeidler}).

\begin{theorem}[Michael Selection Theorem]
\label{MichSel}Let $X$ be a paracompact Hausdorff topological space, and let 
$Y$ be a Banach space. Let $T:X\rightarrow Y$ be a l.s.c. multifunction
having nonempty closed convex images. Then $T$ has a continuous selection.
\end{theorem}

\begin{theorem}[Browder Selection Theorem]
\label{BrSel}Let $X$ be a paracompact Hausdorff topological space, and let $%
Z $ be any topological vector space. Let $T:X\rightarrow Z$ be a
multifunction having nonempty convex images and open preimages. Then $T$ has
a continuous selection.
\end{theorem}

The main goal of this paper is to provide a uniform approach to these four
results. In recent years, these theorems have been generalized in numerous
ways \cite{BC,BD,DT,Ha,HT,Park,PK,Tar,VdV1}. Our main theorems contain as
immediate corollaries a number of these results; however, we do not attempt
to include all of them.

Our approach is based on Michael's \cite{Mich2} and Browder's techniques 
\cite{Br} and the study of abstract convexity structures on topological
spaces originated in works of Michael \cite{Mich3}, Van de Vel \cite{VdV},
Horvath \cite{Ho1}, and others (see Section \ref{Appendix}). Given a
topological (or uniform) space $Y$, Van de Vel introduces the class of
``convex'' sets as a class of subsets of $Y$ closed under intersections.
Horvath defines ``convex hulls'' of finite subsets of $Y$. Michael, on the
other hand, considers an analogue of convex combination functions of vector
spaces: $k_{n}(d_{0},\ldots ,d_{n},a_{0},\ldots
,a_{n})=\sum_{i=0}^{n}d_{i}a_{i},$ where $(d_{0},\ldots ,d_{n})$ is an
element of the $n$-simplex$\ \Delta _{n}$, for certain combinations $%
(a_{0},\ldots ,a_{n})$ of elements of the space. For topological vector
spaces, convex combination functions are continuous with respect to $%
d_{0},\ldots ,d_{n}$. But for locally convex topological vector spaces,
these functions are continuous with respect to $a_{0},\ldots ,a_{n}$ as well.

We follow Michael's approach (Definition \ref{MichDef1})\ because, as we
will see, it is especially convenient for selection and fixed-point
problems. For the former, $Y$ does not have to be ``convex''. Following
Michael, we avoid this situation, while in Van de Vel's and Horvath's
constructions there is the largest ``convex'' set. Our construction extends
Michael's definition in such a way that most of Van de Vel's and Horvath's
fixed-point and selection results are included.

Assuming that $\mathcal{B}$ is a base of the metric uniformity of a space $Y$%
, we state Michael's definition in a form equivalent to the original if $Y$
is compact.

\begin{definition}
\label{MichDef1} A sequence of pairs $\{(M_{n},k_{n})\}$ is a \textit{%
Michael's convex structure} if for all $n\geq 0,\ M_{n}\subset Y^{n+1},\
k_{n}:\Delta _{n}\times M_{n}\rightarrow Y$ ($M_{n}$ can be empty) and the
following is satisfied

(a) if $x\in M_0$, then $k_0(1,x)=x,$

(b) if $x\in M_n,\ n\geq 1,\ i\leq n$, then $\partial _ix\in M_{n-1}$, and
if $t_i=0$ for $t\in \Delta _n$, then $k_n(t,x)=k_{n-1}(\partial
_it,\partial _ix)$, where $\partial _i$ is the operator that omits the $i$th
coordinate,

(c) if $x\in M_n,\ n>i\geq 0$, and $x_i=x_{i+1}$, then for $t\in \Delta _n$%
\begin{equation*}
k_n(t,x)=k_{n-1}(t_1,\ldots ,t_{i-1},t_i+t_{i+1},t_{i+2},\ldots
,t_n,\partial _ix),
\end{equation*}

(d) for each $x\in M_n$, the map $k_n(\cdot ,x)$ is continuous,

(e) for any $U\in \mathcal{B}$ there is a $W\in \mathcal{B}$ such that for
all $n\geq 0,\ t\in \Delta _{n},\ x=(x_{0},\ldots ,x_{n}),\ y=(y_{0},\ldots
,y_{n})\in M_{n}$, we have 
\begin{equation*}
(x_{i},y_{i})\in W,\ 0\leq i\leq n,\Longrightarrow
(k_{n}(t,x),k_{n}(t,y))\in U.
\end{equation*}
\end{definition}

We relax Michael's conditions in several ways (Section \ref{Convexity}).
First, we do not assume that $Y$ is metrizable but only uniform. Second, we
allow the convex combination function to be multivalued. Third, instead of a
sequence of maps $\{k_{n}\}_{n=1}^{\infty }$ connected by conditions (b) and
(c), we use a sole multifunction $C$ from a subset of the set $\Delta (Y)$
of all formal convex combinations of elements of $Y$ into $Y$, which makes
it easier to prove existence of $C$. We consider a certain convexity
multifunction $C$ satisfying continuity conditions (D) and (E) similar to
conditions (d) and (e) above. Condition (E) allow us to carry out most of
the selection and fixed point constructions (Sections \ref{Selection} and 
\ref{FixPoint}) and only as the last step do we consider various continuity
conditions with respect to $t$ (condition (D)), which ensures continuity of
selections and existence of fixed points of continuous maps. The conditions
of Michael's definition do not hold for non-locally convex topological
vector spaces and we have to deal with them in order to obtain Theorems \ref
{Br} and \ref{BrSel}. To resolve this problem, we introduce a second
topology $Z$ on $Y$. As a result, the convexity satisfies the two continuity
requirements, but with respect to two (possibly different) topologies.
Consequently, by choosing an appropriate topological structure on $Y$, we
are able to obtain Theorems \ref{Ka} (for $Z=Y$) and \ref{Br} (for $Y$
discrete) as immediate corollaries of our fixed point result (Theorem \ref
{main}). In the same manner we derive Theorems \ref{MichSel} and \ref{BrSel}
from our selection theorem (Theorem \ref{Sel}).

We shall also use the fact that presence of convex combinations reduces the
question of existence of fixed points for a certain class of multifunctions
on topological spaces to the question of existence of fixed points of
multifunctions on simplexes (see \cite{BD,PK}). This enables us to use the
Brouwer fixed point theorem for $\Delta _{n},$ or, more generally, the
theorem below that is contained in Corollary 2.3 of Gorniewicz \cite{GG1}
(for a stronger result see \cite{Saveliev}). A multivalued map $%
F:X\rightarrow Y$ is called \textit{admissible in the sense of Gorniewicz}
if it is closed valued u.s.c. and there exist a topological space $Z$ and
two single-valued continuous maps $p:Z\rightarrow X,\ q:Z\rightarrow Y$ such
that $p$ is proper and for any $x\in X,$ $(i)\ \ \ p^{-1}(x)$ is acyclic,
and $(ii)\ q(p^{-1}(x))\subset F(x).$ Many fixed point theorems in the
existing literature will be shown to be reducible to this theorem.

\begin{theorem}
\label{Go}Any admissible map in the sense of Gorniewicz and, in particular,
any composition of acyclic multifunctions, $F:\Delta _n\rightarrow \Delta _n$
has a fixed point.
\end{theorem}

Another purpose of this paper is to obtain fixed point theorems for
topological spaces without linear or convex structure. To achieve this goal
we need to show that a given topological space can be equipped with a
convexity structure; that is, to prove existence of convex combination
functions. As a corollary we obtain a generalization of the following
theorem due to Eilenberg and Montgomery \cite{EM}.

\begin{theorem}[Eilenberg-Montgomery Fixed-Point Theorem]
\label{EM}Let $X$ be an acyclic compact ANR, and let $F:X\rightarrow X$ be
an acyclic multifunction. Then $F$ has a fixed point.
\end{theorem}

\section{\-Convexity on Uniform Spaces.\label{Convexity}}

Throughout the paper we fix an infinite cardinal number $\omega $ and an
index set $I$ with $|I|=\omega $. We assume that $\omega $ is large enough
in the sense that $\omega \geq 2^{|X|}$ for each space $X$ involved. Let $%
\Delta _{\omega }$ be the infinite dimensional simplex spanned by the unit
vectors $e_{k},\ k\in I$, of $[0,1]^{\omega }$. For any nonempty subset $K$
of $I$, let $\Delta _{K}$ denote the convex hull of the set $\{e_{k}:k\in
K\} $ in $\Delta _{\omega }$ (i.e. $\Delta _{K}=\{d=(d_{i})_{i\in I}\in
\Delta _{\omega }:i\notin K\Rightarrow d_{i}=0\}$) and let $\Delta _{n}$ be
any $n$-simplex in $\Delta _{\omega }$ spanned by some unit vectors$.$\ We
let $\Delta (Y)=\Delta _{\omega }\times Y^{\omega }$ and assume the
following:

\begin{description}
\item[Convention]  $
\begin{array}{c}
\text{If }(d,a)\in \Delta (Y)\text{ then }d_{i}\in \lbrack 0,1],\ a_{i}\in
Y,\ i\in I,\text{are} \\ 
\text{ the coordinates of }d\in \lbrack 0,1]^{\omega },\ a\in Y^{\omega }.
\end{array}
$
\end{description}

For any $A\subset Y$, we define the set of \textit{all formal convex
combinations} of elements of $A:$%
\begin{equation*}
\Delta (A)=\{(d,a)\in \Delta (Y):d_i\neq 0\Rightarrow a_i\in A,\ i\in I\}.
\end{equation*}

The following are also fixed:

(c1) $Y$ is a uniform space with a minimal (uniform) open base $\mathcal{B}$
(i.e., one with the smallest cardinality), partially ordered by inclusion,

(c2) $Z$ is a topological space on $Y$ (the topology of $Z$ is not
necessarily the uniform topology of $Y$),

(c3) $\mathcal{V}$ is a class of multifunctions that will be specified later,

(c4) $\mathcal{A}\subset 2^{Y}\backslash \{\emptyset \}$ is a class of
subsets of $Y$\textit{\ }called \textit{admissible sets }($\mathcal{A}$ may
be empty)\textit{,}

(c5) $conv:\mathcal{A}\rightarrow 2^{Y}\backslash \{\emptyset \}$ is a
function, and $conv(A)$ is called the \textit{convex hull} of $A\in \mathcal{%
A,}$

(c6) $\mathcal{C}=\{A\in \mathcal{A}:conv(A)\subset A\}\cup \{\emptyset \}$
is the set of \textit{convex subsets} of $Y,$

(c7)$\ Q=\cup_{A\in \mathcal{A}}\Delta (A),$ $Q^{\prime }$ is a subset of 
$\Delta (Y)$ containing $\cup_{A\in \mathcal{A},W\in \mathcal{B}}\Delta
(B(A,W)),$

(c8) $C:Q^{\prime }\rightarrow Y$ is a multifunctions called \textit{convex
combination}.

\begin{definition}[Main Definition]
\label{maindef}The triple $\kappa =(Y,C,Z)$ is called a\textit{\ convexity}
associated with $\mathcal{A,\ }conv\mathcal{,\ V},\ Q^{\prime }$ (this part
will often be suppressed) if the following conditions are satisfied:

(D)\ for any $a\in Y^{\omega }$, if $\Delta _{n}\times \{a\}\subset
Q^{\prime },\ n\geq 0$, then the multifunction $C(\cdot ,a):\Delta
_{n}\rightarrow Z$ belongs to $\mathcal{V}$,

(E)\ for any $U\in \mathcal{B}$, there exists$\ W\in \mathcal{B}$ such that 
\begin{equation*}
C(\Delta (B(A,W)))\subset B(conv(A),U)\text{ for all admissible }A\subset Y.
\end{equation*}
\noindent If $Q^{\prime }=\Delta (Y)$ then the convexity is called \textit{%
global}.
\end{definition}

\begin{proposition}
\label{C'} If $(Y,C,Z)$ is a convexity then the following condition is
satisfied:

$(\gamma )\ \ C(\Delta (A))\subset \overline{conv}(A)$, the closure of $%
conv(A)$ in $Y$, for all admissible $A\subset Y.$
\end{proposition}

\section{A Strong Convexity.\label{Strong}}

Let $F_i:X\rightarrow Y,\ i\in J$, be multifunctions, where $J$ is a
directed set. Then we say that $\{F_i:i\in J\}$ \textit{converges uniformly
on }$N\subset X$ to a multifunction $F:N\rightarrow Y$ if for any $U\in 
\mathcal{B}$, there exists an $i_0\in A$ such that 
\begin{equation*}
F_i(x)\subset B(F(x),U)\text{ for all }x\in N,\ i\in J,\ i>i_0.
\end{equation*}

Let $\Omega $ denote the set of all finite subsets of the index set $I.$ For
a fixed $d\in \Delta _\omega $, we define elements of the product uniformity
of $\{d\}\times Y^\omega $ as follows: for any $W\in \mathcal{B}$, $m\in
\Omega $, we let

\begin{equation*}
\begin{array}{l}
W^{m}=\{((d,a),(d,a^{\prime }))\in \Delta (Y)\times \Delta
(Y):(a_{j},a_{j}^{\prime })\in W,\ j\in m\}, \\ 
B^{\ast }((d,a),W^{m})=\{(d,a^{\prime })\in \Delta (Y):(a_{j},a_{j}^{\prime
})\in W,\ j\in m\}, \\ 
B^{\ast }(S,W^{m})=\cup_{s\in S}B^{\ast }(s,W^{m}),
\end{array}
\end{equation*}
where $(d,a)\in \Delta (Y),\ S\subset \Delta (Y).$

Consider the following conditions on the objects defined in (c1)-(c9) that
loosely correspond to conditions (b)-(e) of Michael's Definition \ref
{MichDef1}:

$(\beta )$ (permutations) if $(d,a),(d^{\prime },a^{\prime })\in Q$ and $%
\sum_{a_i=y}d_i=\sum_{a_i^{\prime }=y}d_i^{\prime }$ for any $y\in Y$, then $%
C(d,a)=C(d^{\prime },a^{\prime }),$

$(\gamma )$ (convex hull) $C(\Delta (A))\subset \overline{conv}(A)$, the
closure in $Y$, for all admissible $A\subset Y,$

$(\delta )=(D)$ ($d$-continuity) for any $a\in Y^{\omega }$, if $\Delta
_{n}\times \{a\}\subset Q^{\prime },\ n\geq 0$, then the multifunction $%
C(\cdot ,a):\Delta _{n}\rightarrow Z$ belongs to $\mathcal{V},$

$(\varepsilon )$ ($a$-continuity) for any $U\in \mathcal{B}$, there exist $%
W\in \mathcal{B},\ m\in \Omega $, such that 
\begin{equation*}
C(B^{\ast }((d,a),W^{m}))\subset B(C(d,a),U)\text{ for all }(d,a)\in Q.
\end{equation*}

\begin{definition}
\label{strdef}The triple $\kappa =(Y,C,Z)$ is called a \textit{strong
convexity} if conditions $(\beta ),\ (\gamma ),\ (\delta )$ and $%
(\varepsilon )$ are satisfied.
\end{definition}

As a direct consequence of the definitions above, we obtain the following
for a strong convexity.

\begin{lemma}
\label{3.2}For any $A\subset Y,\ W\in \mathcal{B},\ m\in \Omega ,$ we have 
\begin{equation*}
\Delta (B(A,W))\subset B^{*}(\Delta (A),W^m).
\end{equation*}
\end{lemma}

\begin{theorem}
\label{C1C2}Conditions $(\alpha ),\ (\gamma )$ and $(\varepsilon )$ imply
condition $(E)$, so any strong convexity is a convexity.
\end{theorem}

\begin{proof}
Let $U\in \mathcal{B}$ be fixed, and let $U^{\prime }\in \mathcal{B}$
satisfy $4U^{\prime }\subset U.$ By $(\varepsilon )$, there exist $W\in 
\mathcal{B},\ m\in \Omega $, such that 
\begin{equation*}
C(B^{\ast }((d,a),W^{m}))\subset B(C(d,a)),U^{\prime })\text{ for all }%
(d,a)\in Q.
\end{equation*}
If $A$ is an admissible set then $\Delta (A)\subset Q$, so this inclusion
holds for all $(d,a)\in \Delta (A).$ Hence 
\begin{equation*}
C(B^{\ast }(\Delta (A),W^{m}))\subset B(C(\Delta (A)),2U^{\prime }).
\end{equation*}
Applying consecutively Lemma \ref{3.2}, the above inclusion and condition $%
(\gamma )$, we obtain condition $(E)$.
\end{proof}

\section{Convexity of Topological Vector Spaces.\label{ContConv}}

\begin{definition}
\label{contconv}Let $\mathcal{V}$ be the class of all single-valued
continuous maps. Then we say that the\textit{\ }convexity is \textit{%
continuous}$.$
\end{definition}

Some examples of spaces with continuous convexity are listed in this and
following sections.

\begin{definition}
A continuous convexity $\kappa =(Y,C,Y)$ (here the topological structures of 
$Y$ and $Z$ coincide) associated with $\mathcal{A}$, $conv$, is called 
\textit{regular} if it is global and for any $y\in Y$, we have $\{y\}\in 
\mathcal{A}$ and $conv(\{y\})=\{y\}$.
\end{definition}

\begin{proposition}
\-\label{4.2}Let $Y$ be a convex subset of a locally convex topological
vector space. Then $(Y,C,Y)$ with $C$ given by 
\begin{equation*}
C(d,a)=\sum_{i}d_{i}a_{i},\quad (d,a)\in \Delta (Y),
\end{equation*}
is a regular (strong) convexity associated with $\mathcal{A}=2^{Y},\
conv(A)=co(A),\ A\subset Y$ (where $co(A)$ is the usual convex hull in a
vector space)$.$
\end{proposition}

If the topological vector space is not locally convex Proposition \ref{4.2}
fails, which motivates the next definition.

\begin{definition}
A continuous convexity $(Y,C,Z)$ is called \textit{discrete} if $Y$ is
discrete (then condition $(E)$ turns into $C(\Delta (A))\subset conv(A)$ for
all $A\in \mathcal{A}$).
\end{definition}

\begin{proposition}
\label{4.3}Let $Z$ be a convex subset of a topological vector space. Then $%
(Y,C,Z)$ ($Y$ is discrete) with $C$ given by 
\begin{equation*}
C(d,a)=\sum_{i}d_{i}a_{i},\quad (d,a)\in \Delta (Y),
\end{equation*}
is a discrete global convexity associated with $\mathcal{A}=2^{Y},\
conv(A)=co(A),\ A\subset Y$.
\end{proposition}

\begin{proof}
Condition $(E)$ is trivially satisfied, because the uniform base $\mathcal{B}
$ of $Y$ consists of only one element $B_{0}=\{(b,b):b\in Y\}$. Now we
observe that $C(\cdot ,a):\Delta _{n}\rightarrow Z$ is continuous as a
linear map on a finite-dimensional space, so $(D)$ of Definition \ref
{maindef} holds for $\mathcal{V}$ the class of all continuous maps.
\end{proof}

\section{\-Constructing a Convexity on a Topological Space.\label{ConvPoints}%
}

Our goal is to construct a regular convexity on a given uniform space.
Throughout this section we assume that the index set $I$ is well ordered.
Let 
\begin{equation*}
\begin{array}{c}
\mathcal{A}=\{\{y\}:y\in Y\}\text{ and} \\ 
Q=\cup_{y\in Y}\Delta (\{y\})=\{(d,a)\in \Delta (Y):\text{for some }y\in
Y,d_{i}\neq 0\Rightarrow a_{i}=y\}
\end{array}
\end{equation*}

\begin{lemma}
\label{7.1}Suppose $Y$ is Hausdorff compact and infinite. Let $\Delta
^{\prime }(Y)$ be the quotient space of $\Delta (Y)$ with respect to the
following equivalence relation: 
\begin{equation*}
(d,a)\sim (d^{\prime },a^{\prime })\text{ if for any }y\in Y,\
\sum_{a_{i}=y}d_{i}=\sum_{a_{i}^{\prime }=y}d_{i}^{\prime }.
\end{equation*}
Let $C:Q\rightarrow Y$ be given by $C(x)=y$ for any $x\in \Delta (\{y\}).$%
\noindent\ Then

(i) $\Delta ^{\prime }(Y)$ is a Hausdorff uniform space and the quotient map 
$p:\Delta (Y)\rightarrow \Delta ^{\prime }(Y)$ is uniformly continuous,

(ii) (a) In $Q$, equivalence classes are $\Delta (\{y\}),\ y\in Y$, (b) $%
p(\Delta (\{y\})=y,\ y\in Y$, and (c) $p(Q)=Y$.
\end{lemma}

\begin{proof}
First, $\Delta (Y)$ is normal, so $\Delta ^{\prime }(Y)$ is Hausdorff. Next,
for any $(d,a)\in \Delta (Y),\ y\in Y$, let $S_{y}(d,a)=\sum_{a_{i}=y}d_{i}.$
Then the set 
\begin{equation*}
R=\{((d,a),(d^{\prime },a^{\prime })):\text{for any }y\in Y,\ \
S_{y}(d,a)=S_{y}(d^{\prime },a^{\prime })\}\text{ }
\end{equation*}
determines the equivalence relation $\sim .$ To ensure that $\Delta ^{\prime
}(Y)$ has a uniformity, we need to check that for any $W\in \mathcal{B,}$
any $\varepsilon >0,\ m\in \Omega ,$ there are $V\in \mathcal{B},\ \delta
>0,\ n\in \Omega ,$ such that 
\begin{equation*}
V^{\delta ,n}+R+V^{\delta ,n}\subset R+W^{\varepsilon ,m}+R,
\end{equation*}
which means that the equivalence relation $\sim $ and the uniformity of $%
\Delta (Y)$ are weakly compatible (see \cite[p. 24]{Jam}), so (i) holds. We
will show that $\Delta (Y)\times \Delta (Y)\subset R+W^{\varepsilon ,m}+R.$
Consider 
\begin{equation*}
\begin{array}[t]{ll}
R+W^{\varepsilon ,m}+R & =\{((d,a),(d^{\prime },a^{\prime })):\text{there is
an }(r,b)\in \Delta (Y)\text{ with } \\ 
& \qquad ((d,a),(r,b))\in R,\ ((r,b),(d^{\prime },a^{\prime }))\in
W^{\varepsilon ,m}+R\} \\ 
& =\{((d,a),(d^{\prime },a^{\prime })): \\ 
& \text{ \ \ }\exists \text{ }(r,b)\in \Delta (Y)\text{ with }%
S_{y}(d,a)=S_{y}(r,b),\forall y\in Y, \\ 
& \text{ \ \ }\exists \text{ }(r^{\prime },b^{\prime })\in \Delta (Y)\text{
with }S_{y}(d^{\prime },a^{\prime })=S_{y}(r^{\prime },b^{\prime }),\forall
y\in Y, \\ 
& \text{ \ \ and }\mid r_{i}-r_{i}^{\prime }\mid <\varepsilon ,\
(b_{i},b_{i}^{\prime })\in W,\ i\in m\} \\ 
& =\{((d,a),(d^{\prime },a^{\prime })):\text{there are }(r,b),(r^{\prime
},b^{\prime })\in \Delta (Y)\text{ with} \\ 
& \text{ \ \ (1) }\mid r_{i}-r_{i}^{\prime }\mid <\varepsilon ,\
(b_{i},b_{i}^{\prime })\in W,\ i\in m, \\ 
& \text{ \ \ (2) }S_{y}(d,a)=S_{y}(r,b),\forall y\in Y, \\ 
& \text{ \ \ (3) }S_{y}(d^{\prime },a^{\prime })=S_{y}(r^{\prime },b^{\prime
}),\forall y\in Y\}.
\end{array}
\end{equation*}
Suppose $((d,a),(d^{\prime },a^{\prime }))\in \Delta (Y)\times \Delta (Y)$.
Assume for simplicity that $m=\{0,1,...,M\}$. Choose $z\in Y$ such that $%
d_{i}\neq 0\Rightarrow z\neq a_{i}$ and $d_{i}^{\prime }\neq 0\Rightarrow
z\neq a_{i}^{\prime }$ for all $i\in I.$ To get $(r,b),(r^{\prime
},b^{\prime })\in \Delta (Y)$ as above, let 
\begin{equation*}
r_{i}=r_{i}^{\prime }=0,\ b_{i}=b_{i}^{\prime }=z,\qquad i\in m=\{0,1,\ldots
,M\},
\end{equation*}
Then condition (1) above is satisfied. Now, let 
\begin{equation*}
r_{M+i+1}=d_{i},\ r_{M+i+1}^{\prime }=d_{i}^{\prime },\ b_{M+i+1}=a_{i},\
b_{M+i+1}^{\prime }=a_{i}^{\prime },\ i=0,1,\ldots
\end{equation*}
(it means that we obtain $(r,b)$ and $(r^{\prime },b^{\prime })$ by
``shifting'' coordinates of $(d,a)$ and $(d^{\prime },a^{\prime })$ $M$
steps to the right). Then conditions (2) and (3) above are also satisfied.
Therefore we have 
\begin{equation*}
R+W^{\varepsilon ,m}+R=\Delta (Y)\times \Delta (Y).
\end{equation*}

If $(d,a)\in Q$ and $d_{i}\neq 0$ then $a_{i}=y$ for all $i$ and some $y\in
Y $, and $p(d,a)=y.$ Then for any $V\in \mathcal{B},\ n\in \Omega ,\
\varepsilon >0,$ we have

\begin{equation*}
\begin{array}{ll}
B^{\ast }((d,a),V^{\varepsilon ,n})\cap Q & =\{(d^{\prime },a^{\prime })\in
\Delta (Y):\text{for some }\ y^{\prime }\in B(y,V), \\ 
& \qquad \qquad d_{i}^{\prime }\neq 0\Rightarrow a_{i}^{\prime }=y^{\prime },%
\text{ and }\mid d_{i}-d_{i}^{\prime }\mid <\varepsilon ,\ i\in n\} \\ 
& =\cup_{y^{\prime }\in B(y,V)}\{(d^{\prime },a^{\prime })\in \Delta
(Y):d_{i}^{\prime }\neq 0\Rightarrow a_{i}^{\prime }=y^{\prime }, \\ 
& \qquad \qquad \qquad \qquad \qquad \mid d_{i}-d_{i}^{\prime }\mid
<\varepsilon ,\ i\in n\}.
\end{array}
\end{equation*}
Therefore we have 
\begin{equation*}
p(B^{\ast }((d,a),V^{\varepsilon ,n})\cap Q)=B(y,V),
\end{equation*}
and we conclude that $p(Q)=Y$.
\end{proof}

\begin{theorem}
\label{7.4}Let $Y$ be a Hausdorff compact ANR and let $\mathcal{V}$ be the
set of all u.s.c. multifunctions with values either singletons or $Y$. Then
there exists a convexity $\kappa =(Y,C,Y)$ associated with $\mathcal{A}%
=\{\{y\}:y\in Y\},\ conv(\{y\})=\{y\}$ for all $y\in Y$.
\end{theorem}

\begin{proof}
By Lemma \ref{7.1} $\Delta ^{\prime }(Y)$ is Hausdorff and $p(Q)=Y$ is
compact. Hence $p(Q)$ is closed in $\Delta ^{\prime }(Y).$ Then, there
exists an open neighborhood $N_{V}\subset \Delta ^{\prime }(Y)$ of $p(Q)$
and a continuous function $C:\overline{N_{V}}\rightarrow Y$ that extends $%
Id_{Y}$ \ and $C$ is uniformly continuous. Now if we extend $C$ on the whole 
$\Delta ^{\prime }(Y)$ by putting $C(d,a)=Y$ for $(d,a)\notin N_{V},$ then $%
C $ is uniformly u.s.c. (since $N_{V\text{ }}$ is open). It is routine to
check that $C$ satisfies conditions $(\beta )$ -- $(\delta )$.
\end{proof}

\section{\-Preliminaries From General Topology.\label{GenTop}}

Presenting necessary definitions we mostly follow Engelking \cite{En}.

\begin{definition}
The \textit{weight} $w(Y)$ is the cardinality of $\mathcal{B}$. Let $\varphi
(Y)$ be the largest cardinal number $\mu \leq \omega $ such that the
intersection of any family of elements of $\mathcal{B}$ whose cardinality is
less then $\mu $ contains an element of $\mathcal{B}$.
\end{definition}

Note that if $\mathcal{B=}\{U_\beta :\beta <\mu \}$, where $\mu $ is an
infinite ordinal, is ordered by inclusion (in this case $Y$ is called a $\mu 
$\textit{-metrizable} space \cite{Hr}), then $\varphi (Y)=w(Y)^{+},$ where $%
\alpha ^{+}$ stands for the least cardinal number larger that $\alpha .$

\begin{definition}
For a topological space $X$, the \textit{Lindel\"of number} $l(X)$ is the
least cardinal number $\lambda $ such that every open cover of $X$ has a
subcover whose cardinality is less than $\lambda $ (''at most'' in \cite{En}%
). Let $l^{\prime }(X)$ be the largest cardinal number $\mu \leq \omega $
such that any open cover of $X$ whose cardinality is less than $\mu $ has a
finite subcover. Let $p(X)$ be the largest cardinal number $\kappa \leq
\omega $ such that any open cover of $X$ whose cardinality is less than $%
\kappa $ has a locally finite open refinement.
\end{definition}

Then $X$ is known \cite{St} as \textit{finally }$\lambda $\textit{-compact}
(or $\lambda $-\textit{Lindel\"{o}f }\cite{Hr}) and \textit{initially }$\mu $%
\textit{-compact}, respectively.

\begin{definition}
For the uniform space $Y$, let $l_u(Y)$ be the least cardinal number $%
\lambda $ such that for every $V\in \mathcal{B}$, the cover $\{B(x,V):x\in
Y\}$ has a subcover whose cardinality is less than $\lambda .$
\end{definition}

The proof of the following is similar to \cite[Theorem 3.1.23]{En}.

\begin{proposition}
\label{corgen}Suppose $\mathcal{B}=\{U_{i}:i\in J\}$ is partially ordered by
inclusion and $\varphi (Y)\geq l(Y)$. Then any net $\{y_{i}:i\in J\}$ in $Y$
has a convergent subnet.
\end{proposition}

\section{\-The Main Selection Theorems.\label{Selection}}

Since Theorems \ref{Ka} - \ref{BrSel} deal with two types of multifunctions
(u.s.c. and l.s.c.), we shall consider maps $T:X\rightarrow Y$ and $%
R:Z\rightarrow X$ of either kind between the two spaces ($Z=Y$) and a fixed
point of their composition $T\circ R:Y\rightarrow Y.$ The diagram below
illustrates proofs in Sections \ref{Selection} and \ref{FixPoint}.

If a multifunction $T:X\rightarrow Y$ has admissible images then its \textit{%
convex hull} $conv(T):X\rightarrow Y$ is given by 
\begin{equation*}
conv(T)(x)=conv(T(x)),\quad x\in X.
\end{equation*}

Michael proves his selection theorem for his convex structures \cite{Mich3}
by considering a sequence of ``almost continuous'' selections, while for
locally convex topological vector spaces he constructs \cite{Mich2} a
sequence of continuous ``almost selections''. The former yields sharper
selection results, the latter requires an additional restriction on the
convexity (some neighborhood of an admissible set is admissible, as in
Proposition \ref{MichConv}), but allows us to proceed directly to fixed
point theorems.

\begin{theorem}[Almost Selection Theorem]
\label{selection}Let $X$ be a normal topological space, $\kappa =(Y,C,Z)$ a
convexity, $Y^{\prime }$ a subset of $Y$. Suppose also that:

(i) $T:X\rightarrow Y$ is admissible-valued l.s.c. and $T(x)\cap Y^{\prime
}\neq \emptyset $ for any $x\in X$,

(ii) $p(X)\geq l_u(Y^{\prime }).$

\noindent Then for any $U\in \mathcal{B,}$ there exist $a\in (Y^{\prime
})^{\omega }$ and a continuous function $f:X\rightarrow \Delta _{\omega }$
satisfying the following conditions:

(a) for any $x\in X$, there is an open neighborhood $G$ of $x$ such that $%
f(G)\subset \Delta _n\subset \Delta _\omega $ for some $n\geq 0$,

(b) $f(X)\times \{a\}\subset Q^{\prime },$ and

(c) $C(f(x),a)\subset B(conv(T)(x),U)$ for all $x\in X.$

\noindent If, in addition, $Y^{\prime }$ is admissible, then 
\begin{equation*}
C(f(x),a)\subset \overline{conv}(Y^{\prime }).
\end{equation*}

If, moreover, $\ $(ii$^{\prime }$) $l^{\prime }(X)\geq l_{u}(Y^{\prime }),$
\noindent then we have (a$^{\prime }$) $f(X)\subset \Delta _{n}\subset
\Delta _{\omega }$ for some $n\geq 0$.
\end{theorem}

\begin{proof}
Let $U\in \mathcal{B.}$ Then condition $(E)$ reads as follows: there exists $%
W\in \mathcal{B}$ such that 
\begin{equation}
C(\Delta (B(A,W)))\subset B(conv(A),U)\text{ for all admissible }A\subset Y.
\label{1}
\end{equation}
Let $M=\{B(y,W):y\in Y^{\prime }\}$. By definition of $l_{u}(Y^{\prime })$, $%
M$ has a subcover $M^{\prime }$ with $|M^{\prime }|<l_{u}(Y^{\prime })$. But 
$T:X\rightarrow Y$ is l.s.c., so $N=\{T^{-1}(G):G\in M^{\prime }\}$ consists
of open sets. Moreover, if $x\in X$ then $T(x)\cap G\neq \emptyset $ for
some $G\in M^{\prime }$. Hence $x$ belongs to $T^{-1}(G)$, so $N$ is an open
cover of $X$. By (ii), we have $|N|=|M^{\prime }|<l_{u}(Y^{\prime })\leq
p(X) $. Therefore, by definition of $p(X)$, $N$ has a locally finite open
refinement $N^{\prime }.$ Since $|I|=\omega >2^{|X|}$, we can assume that $%
N^{\prime }=\{Q_{k}:k\in I\}$ ($Q_{k}=\emptyset $ for some $k\in I$). From
the definitions of $M,\ N,\ N^{\prime },$ it follows that for all $k\in I,$ $%
Q_{k}\subset T^{-1}(G_{k}),$ where $G_{k}=B(a_{k},W)$ for some $a_{k}\in
Y^{\prime }$ (here we assign an index to $G_{k}$ and $a_{k}$ according to
this inclusion). Then let $a=(a_{i})_{i\in I}\in (Y^{\prime })^{\omega }.$
From the fact that $X$ is normal and Michael's Lemma it follows that there
exists a partition of unity subordinate to $N^{\prime }$, i.e., there are
continuous functions $f_{k}:X\rightarrow \lbrack 0,1],\ k\in I$, satisfying 
\begin{equation*}
\begin{array}{c}
f_{k}(x)=0\text{ for any }x\notin Q_{k},\ k\in I\text{, and } \\ 
\sum_{k\in I}f_{k}(x)=1\text{ for any }x\in X.
\end{array}
\end{equation*}
Now let 
\begin{equation*}
f(x)=\sum_{k\in I}f_{k}(x)e_{k},
\end{equation*}
where $e_{k},\ k\in I$, are the vertices of $\Delta _{\omega }.$ Then $%
f:X\rightarrow \Delta _{\omega }$ is a continuous function. Since $N^{\prime
}$ is locally finite, for each $x\in X,$ there are a neighborhood $G$ of $x$
and a finite set $S\subset I$ such that $f_{k}|_{G}$ is nonzero only for $%
k\in S$. Therefore we have $f(G)\subset \Delta _{S}$, so (a) is satisfied.
Moreover, if (ii$^{\prime }$) holds, then $N^{\prime }$ is finite and so (a$%
^{\prime }$) holds$.$

Take $x\in X$. Let 
\begin{equation*}
K=\{k\in I:T(x)\cap G_{k}\neq \emptyset \}\text{ and }A_{K}=\{a_{k}:k\in
K\}\subset Y^{\prime }.
\end{equation*}
If $k\in I$ is such that $f_{k}(x)\neq 0$ then $x\in Q_{k}\subset
T^{-1}(G_{k})$, or $T(x)\cap G_{k}\neq \emptyset $. Hence $k\in K$. By
definition of $f$, this implies that 
\begin{equation}
f(x)=\sum_{k\in I}f_{k}(x)e_{k}=\sum_{k\in K}f_{k}(x)e_{k}\in \Delta _{K},
\label{2}
\end{equation}
where $e_{k},k\in I$, are the unit vectors of $[0,1]^{\omega }$. Next,
consider 
\begin{equation}
\Delta _{K}\times \{a\}=\{(d,a):d_{i}\neq 0\Rightarrow a_{i}\in
A_{K}\}\subset \Delta (A_{K}).  \label{3}
\end{equation}
By definition of $K,\ $ we have $a_{k}\in B(T(x),W)$ for all $k\in K$, or 
\begin{equation}
A_{K}\subset B(T(x),W)\cap Y^{\prime }.  \label{4}
\end{equation}
From (\ref{2})-(\ref{4}), it follows that 
\begin{equation}
f(x)\times \{a\}\in \Delta _{K}\times \{a\}\subset \Delta (A_{K})\subset
\Delta (B(T(x),W)\cap Y^{\prime }).  \label{5}
\end{equation}
Since $T(x)$ is admissible, $\Delta (B(T(x),W))\subset Q^{\prime }$. Hence
by (\ref{5}), we have $f(X)\times \{a\}\subset Q^{\prime }$ and, therefore, $%
C(f(x),a)$ is well defined. Thus, from (\ref{5}) and (\ref{1}) we have (c).

To finish the proof, we notice that (\ref{5}) implies that $f(x)\times
\{a\}\in \Delta (Y^{\prime })$. Therefore we have 
\begin{equation*}
C(f(x),a)\subset C(\Delta (Y^{\prime }))\subset \overline{conv}(Y^{\prime }).
\end{equation*}
\end{proof}

\begin{corollary}[Continuous Almost Selection Theorem]
\label{ApprSel}Let $X$ be a normal topological space, $\kappa =(Y,C,Z)$ a
continuous convexity, $Y^{\prime }$ a subset of $Y.$ Suppose also that

(i) $T:X\rightarrow Y$ is l.s.c. with admissible images and $T(x)\cap
Y^{\prime }\neq \emptyset $ for any $x\in X$,

(ii) $p(X)\geq l_{u}(Y^{\prime }).$

\noindent Then for any $U\in \mathcal{B}$, there exists a continuous $V$%
-almost selection for the multifunction $conv(T):X\rightarrow Z$, i.e.,
there is a continuous function $g:X\rightarrow Z$ such that 
\begin{equation*}
g(x)\in B(conv(T(x)),U)\text{ for all }x\in X.
\end{equation*}
If, moreover, $Y^{\prime }$ is admissible, then we have 
\begin{equation*}
g(X)\subset \overline{conv}(Y^{\prime }).
\end{equation*}
\end{corollary}

\begin{definition}
\label{convbase}We say that the convexity $\kappa =(Y,C,Z)$ has a \textit{%
convex uniform base} $\mathcal{B}$ if 
\begin{equation*}
y\in Y,\ U\in \mathcal{B},\ D\in \mathcal{C}\Rightarrow B(y,U)\cap D\in 
\mathcal{C}
\end{equation*}
\end{definition}

\begin{theorem}[Continuous Selection Theorem]
\label{Sel}Let $X$ be a normal topological space, $Y$ be a complete uniform
space, $\kappa =(Y,C,Z)$ a continuous convexity with a countable convex
uniform base $\mathcal{B}$, and suppose that the uniform topology of $Y$ is
finer than the topology of $Z$. Suppose also that

(i) $T:X\rightarrow Y$ is l.s.c. with nonempty convex images,

(ii) $p(X)\geq l_u(Y).$

\noindent Then the multifunction $\overline{T}:X\rightarrow Z$ has a
continuous selection, i.e., there is a continuous map $g:X\rightarrow Z$
such that 
\begin{equation*}
g(x)\in \overline{T(x)}\text{ for all }x\in X\text{ (closure in }Y\text{)}.
\end{equation*}
\end{theorem}

\begin{proof}
Let $\mathcal{B}=\{U_{1},U_{2},\ldots \}$. Without loss of generality we can
assume that 
\begin{equation}
2U_{n+1}\subset U_{n},\ n=1,2,\ldots .  \label{m4}
\end{equation}
Then according to Corollary \ref{ApprSel} (with $Y^{\prime }=Y$), for any
nonempty convex-valued l.s.c. map $G:X\rightarrow Y,$ for any $U\in \mathcal{%
B}$, there is a continuous $g:X\rightarrow Z$ with 
\begin{equation}
g(x)\in B(G(x),U)\text{ for all }x\in X.  \label{m5}
\end{equation}
We inductively construct a sequence of continuous functions $%
g_{n}:X\rightarrow Z,\ n=1,2,\ldots $, such that 
\begin{equation}
g_{n}(x)\in B(T(x),U_{n+1})\text{ for all }x\in X,\ n=1,2,\ldots ,
\label{m6}
\end{equation}
\begin{equation}
g_{n}(x)\in B(g_{n-1}(x),U_{n-1})\text{ for all }x\in X,\ n=2,3,\ldots .
\label{m7}
\end{equation}
By (\ref{m5}), there is a$\ g_{1}$ so that (\ref{m6}) holds for $n=1.$
Assume that we have constructed $g_{1},\ldots ,g_{n-1}$ satisfying these
conditions. Then let 
\begin{equation}
G(x)=B(g_{n-1}(x),U_{n})\cap T(x).  \label{m10}
\end{equation}
Then $G$ is l.s.c.. By (\ref{m6}) for $n-1$, we have $g_{n-1}(x)\in
B(T(x),U_{n})$, so $G(x)$ is nonempty, and it is convex because the base in
convex. Therefore by (\ref{m5}), there is a continuous map $%
g_{n}:X\rightarrow Z$ with 
\begin{equation}
g_{n}(x)\in B(G(x),U_{n+1})\text{ for all }x\in X.  \label{m12}
\end{equation}
Then there is a $y\in G(x)$ such that 
\begin{equation}
(g_{n}(x),y)\in U_{n+1}.  \label{m8}
\end{equation}
By (\ref{m10}), we have $G(x)\subset T(x),$ so (\ref{m12}) implies (\ref{m6}%
). From (\ref{m10}), it also follows that $G(x)\subset B(g_{n-1}(x),U_{n})$.
Therefore we have $(y,g_{n-1}(x))\in U_{n},$ and now, from (\ref{m8}), (\ref
{m11}), and (\ref{m4}), it follows that $(g_{n}(x),g_{n-1}(x))\in U_{n-1}.$
Hence (\ref{m7}) holds. Thus, we have constructed a sequence $%
\{g_{n}:n=1,2,\ldots \}$ satisfying the required conditions.

Now (\ref{m4}) implies that this is a Cauchy sequence. Therefore $%
\{g_{n}:n=1,2,\ldots \}$ converges to a map $g:X\rightarrow Z.$ And from (%
\ref{m6}), it follows that $g(x)\in \overline{T(x)}$ for all $x\in X.$ To
finish the proof we observe that $g:X\rightarrow Z$ is the uniform limit of $%
\{g_{n}\}$ with respect to the uniformity of $Y$. Therefore $g$ is
continuous, since $Y$ is finer than $Z$.
\end{proof}

\section{More Selection Theorems.\label{MoreSel}}

\begin{theorem}[Michael-Type Selection Theorem]
\label{KakTypeSel}Let $X$ be a (Hausdorff) paracompact space, $\kappa
=(Y,C,Y)$ a continuous convexity, $Y$ complete with a countable convex
uniform base, $T:X\rightarrow Y$ l.s.c. with nonempty convex images. Then $%
\overline{T}$ has a continuous selection$.$
\end{theorem}

\begin{proof}
In Theorem \ref{Sel}, we let $Z=Y$ and notice that $p(X)=\omega \geq
l_{u}(Y) $, so condition (ii) of the theorem holds$.$
\end{proof}

This corollary implies the following results: (1) The Michael Selection
Theorem \ref{MichSel} for Banach spaces, (2) Theorem 1.3 of Michael \cite
{Mich3} for Michael's convex structures (with two additional assumptions:
(a) for any $x\in X,$ any $W\in \mathcal{B,}$ $B(\varphi (x),W)$ is
M-admissible, (b) for any $y\in E,$ any $W\in \mathcal{B,}$ $B(y,W)$ is
M-convex), (3) Theorem 3.3 of Horvath \cite{Ho1} for H-spaces, (4) Part (2)
of Theorem 3.5 of Van de Vel \cite[p. 440]{VdV} (or part (b) of Theorem 4.3
of \cite{VdV1}). (Here (3) and (4) may be looked at as selection theorems
for a l.s.c. map with convex range.)

An example of a pair $X,Y$ that satisfies conditions (ii) of Theorem \ref
{Sel} but is not covered by the Michael Selection Theorem: $X$ is normal but
not necessarily paracompact and $Y$ is precompact. Another example: $X$ is
countably paracompact, or countably compact, such as the space $W$ of all
countable ordinal numbers, and $Y$ is a separable Banach space, such as $l_2$
or $C[0,1].$ Nedev \cite{Nedev} proved a Michael-type selection theorem for $%
X=W$ and $Y$ a reflexive Banach space.

\begin{theorem}[Browder-Type Selection Theorem]
\label{BrTypeSel}Let $X$ be a (Hausdorff) paracompact space, and suppose $%
(Y,C,Z)$ is a discrete convexity, and $T:X\rightarrow Y$ has nonempty
admissible images and open fibers. Then $conv(T):X\rightarrow Z$ has a
continuous selection.
\end{theorem}

\begin{proof}
In Corollary \ref{ApprSel} we let $Y$ be discrete and $Y^{\prime }=Y$. Then $%
T:X\rightarrow Y$ is l.s.c., so condition (i) of the theorem holds. But
since $p(X)=\omega \geq l_{u}(Y)$, (ii) also holds. Finally, we notice that
an almost selection with respect to the discrete topology is in fact a
selection.
\end{proof}

This corollary implies the following results: (1) The Browder Selection
Theorem \ref{BrSel} for topological vector spaces, (2) Theorem 3.2 of
Horvath \cite{Ho1} for H-spaces, (3) Van de Vel's version of the Browder
Selection Theorem \cite[IV.3.23, p. 450]{VdV}.

\section{Classes of Maps with Fixed Point Conditions. \label{F(Y,X)}}

Motivated by Ben-El-Mechaiekh and Deguire \cite{BD}, Park and Kim \cite{PK}
define an abstract class $\mathcal{U}_{c}^{\kappa }(Y,X)$ of maps that helps
reduce the fixed point problem to the one for multifunctions on simplexes.
\noindent Examples of $\mathcal{U}_{c}^{\kappa }(Y,X)$ are, for instance,
the classes of all u.s.c. multifunctions with compact convex values in
locally convex topological vector spaces, acyclic maps, and approachable
maps. Taking this one step further, we introduce the following.

\begin{definition}
\label{FYX}Let $X$ be a topological space, $\kappa =(Y,C,Z)$ be a convexity
and $Y^{\prime }$ a subset of $Y$ (or $Z$). Then the class $\mathcal{F}%
_{\kappa }(Y^{\prime },X)$ is defined as the class of all multifunctions $%
F:Z\supset Y^{\prime }\rightarrow X$ such that for any simplex $\Delta
_{n}\subset \Delta _{\omega }$, any $a\in (Y^{\prime })^{\omega }$, and any
continuous function $f:X\rightarrow \Delta _{n}$, the multifunction 
\begin{equation*}
f\circ F\circ C\mid _{\Delta _{n}\times \{a\}}:\Delta _{n}\rightarrow \Delta
_{n}
\end{equation*}
has a fixed point.
\end{definition}

By Definitions \ref{maindef} and \ref{contconv}, $C(\cdot ,a):\Delta
_{n}\rightarrow Z,\ V\in B$, is continuous when $\kappa =(Y,C,Z)$ is a
continuous convexity. Therefore the Brouwer Fixed Point Theorem implies:

\begin{proposition}
\label{6.2}If $\kappa =(Y,C,Z)$ is a continuous global convexity then $%
\mathcal{F}_{\kappa }(Y,X)$ contains all continuous functions $\psi
:Z\rightarrow X.$
\end{proposition}

If $\mathcal{V}$ is the set of all acyclic maps, then we say that a\textit{\ 
}convexity associated with $\mathcal{V}$ is \textit{acyclic}$.$ Then Theorem 
\ref{Go} implies:

\begin{proposition}
\label{6.3}If $\kappa =(Y,C,Z)$ is a global acyclic convexity then $\mathcal{%
F}_{\kappa }(Y,X)$ contains all admissible in the sense of Gorniewicz (and,
therefore, all acyclic) multifunctions $F:Z\rightarrow X$.
\end{proposition}

\begin{proposition}
\label{Y'}If $\kappa =(Y,C,Y)$ is a global acyclic convexity and $Y^{\prime
} $ is a closed convex subset of $Y,$ then $\mathcal{F}_{\kappa }(Y^{\prime
},X)$ contains all u.s.c. multifunctions $F:Y^{\prime }\rightarrow X$ such
that $F(x)$ is compact and acyclic for all $x\in Y^{\prime }.$
\end{proposition}

\begin{proposition}
\label{Uc}If $\kappa =(Y,C,Y)$ is a global continuous convexity then we have 
\begin{equation*}
\mathcal{U}_{c}^{\kappa }(Y,X)\subset \mathcal{F}_{\kappa }(Y,X).
\end{equation*}
\end{proposition}

\section{The Main Fixed Point Theorems.\label{FixPoint}}

\begin{theorem}[Almost Fixed Point Theorem]
\label{apprmain} Let $X$ be a normal topological space, $\kappa =(Y,C,Z)$ a
global convexity. Suppose also that the following holds:

(i) $R\in \mathcal{F}_{\kappa }(Y,X),$

(ii) $T:X\rightarrow Y$ is admissible-valued l.s.c.$,$

(iii) $l^{\prime }(X)\geq l_u(Y).$

\noindent Then the multifunction $conv(T)\circ R:Y\rightarrow Y$ has a $U$-%
\textit{almost fixed point} for any $U\in \mathcal{B}$, i.e., there is a 
\begin{equation*}
y\in B(conv(T)\circ R(y),U).
\end{equation*}
\end{theorem}

\begin{proof}
By Theorem \ref{selection} with (ii$^{\prime }$), for any $U\in \mathcal{B}$
there exist $a\in Y^{\omega }$ and a continuous function $f:X\rightarrow
\Delta _{J}$ with $J$ finite such that $C(f(x),a)\subset B(conv(T)(x),U)$
for all $x\in X.$ Since the convexity is global, we can define a map $\Psi
:\Delta _{J}\rightarrow Y$ by $\Psi =C|_{\Delta _{J}\times \{a\}}.$ Since $%
R\in \mathcal{F}_{\kappa }(Y,X)$, the multifunction $F=f\circ R\circ \Psi
:\Delta _{J}\rightarrow \Delta _{J}$ has a fixed point, that is, there
exists a $d\in f(R(\Psi (d)))$. Therefore, there exist such $x\in X$ and $%
y\in Y$ that $y\in \Psi (f(x))$ and $x\in R(y).$ Hence 
\begin{equation*}
y\in C(f(x),a)\subset B(conv(T)(x),U)\subset B(conv(T)\circ R(y),U).
\end{equation*}
\end{proof}

\begin{theorem}[Fixed Point Theorem]
\label{main}Let $X$ be a normal topological space, $\kappa =(Y,C,Z)$ a
global convexity. Suppose that the following conditions are satisfied:

(i) $R\in \mathcal{F}_\kappa (Y,X),$

(ii) $T:X\rightarrow Y$ is admissible-valued l.s.c.$,$

(iii) $S:Y\rightarrow Y$ is closed-valued u.s.c. and $conv(T)\circ R\subset
S $,

(iv) $\varphi (Y)\geq l(Y),\ l^{\prime }(X)\geq l_u(Y).$

\noindent Then $S$ has a fixed point.
\end{theorem}

\begin{proof}
Suppose $\mathcal{B}=\{U_{i}:i\in J\}.$ By Theorem \ref{apprmain}, for any $%
i\in J$, there is a $y_{i}\in B(conv(T)\circ R(y_{i}),U_{i}).$ Now we need
to show that $y^{\ast }\in S(y^{\ast })$ for some $y^{\ast }\in Y$. By
Corollary \ref{corgen}, the net $\{y_{i}:i\in J\}$ has a convergent subnet.
Therefore, we can simply assume that $y_{i}\rightarrow y^{\ast }\in Y$. Then
for any $i\in J$, there exists a $k(i)\in J$ such that $y_{k}\in B(y^{\ast
},U_{i})$ for all $k>k(i)$, or 
\begin{equation}
y^{\ast }\in B(y_{k},U_{i})\text{ for all }k>k(i).  \label{7}
\end{equation}
By (iii), $y_{k}\in B(S(y_{k}),U_{k})$ for all $k\in J$. Therefore 
\begin{equation}
y_{k}\in B(S(y_{k}),U_{i})\text{ for all }k>i.  \label{8}
\end{equation}
Since $y_{i}\rightarrow y^{\ast }$ and the multifunction $S$ is u.s.c., then
for any $i\in J$, there exists a $j(i)\in J$ such that 
\begin{equation}
S(y_{k})\subset B(S(y^{\ast }),U_{i})\text{ for all }k>j(i).  \label{9}
\end{equation}
Now using (\ref{7}), (\ref{8}), (\ref{9}), we obtain 
\begin{equation*}
y^{\ast }\in B(S(y^{\ast }),3U_{i})\text{ for all }i\in J.
\end{equation*}
Hence $y^{\ast }\in \overline{S(y^{\ast })}=S(y^{\ast })$, and the proof is
complete.
\end{proof}

\section{More Fixed Point Theorems. \label{Kakutani}}

One can obtain Kakutani and Browder type theorems for spaces with
generalized convexity by repeating arguments that work for topological
vector spaces, see \cite{Ho1,Ho2,Park,Pas,Tar}. But then the results apply
only to convex-valued multifunctions. Therefore they are never stronger than
the Eilenberg-Montgomery Theorem \ref{EM}, which deals with acyclic-valued
multifunctions. So we apply a version of the Eilenberg-Montgomery Theorem,
Theorem \ref{Go}, or use $\mathcal{F}_{\kappa }(X,X)$ to obtain sharper
results.

\begin{theorem}[Kakutani-Type Fixed Point Theorem]
\label{corKa}Let $X$ be a (Hausdorff) compact topological space, $\kappa
=(Y,C,Y)$ a regular convexity, $f:X\rightarrow Y$ a (single-valued)
continuous map and $R\in \mathcal{F}_{\kappa }(Y,X)$ u.s.c.. Then $\overline{%
f\circ R}$ has a fixed point.
\end{theorem}

\begin{proof}
In Theorem \ref{main} we let $T=f,\ S=\overline{f\circ R}.$ Then (i) and
(ii) hold because the convexity is regular: $conv(\{y\})=\{y\}$. Also, since 
$Y$ is compact, we have $l^{\prime }(X)\geq \aleph _{0}=l(Y)$ and $\varphi
(Y)\geq \aleph _{0}=l(Y)$, so (iv) holds.
\end{proof}

This corollary implies the following results: (1) The Kakutani Theorem \ref
{Ka} for locally convex topological vector spaces, (2) Corollary to Theorem
6 of Horvath \cite{Ho2} for l.c.-spaces, (3) Theorem 6.15 of Van de Vel 
\cite[p. 498]{VdV}, (4) Theorem 1 of Had\v{z}i\'{c} \cite{Ha} for H-spaces
(with the assumption that all points are H-convex and that space is normal).

\begin{theorem}[Browder-Type Fixed Point Theorem]
\label{corBr}Let $X$ be compact, $(Y,C,Z)$ be a discrete global convexity$.$
Suppose that

(i) $R\in \mathcal{F}_\kappa (Y,X),$

(ii) $T:X\rightarrow Y$ is a multifunction with admissible images and open
fibers.

\noindent Then the multifunction $conv(T)\circ R$ has a fixed point.
\end{theorem}

\begin{proof}
Notice that if $Y$ is a discrete uniform space, then a multifunction $%
T:X\rightarrow Y$ has open fibers if and only if it is l.s.c.. Therefore
conditions (i), (ii) of this theorem are exactly conditions (i), (ii) of
Theorem \ref{main}. Condition (iii) follows from the fact that $%
S=conv(T)\circ R:Y\rightarrow Y$ is u.s.c. with respect to the discrete
topology. Also, since $X$ is compact and $Y$ is discrete, we have $l^{\prime
}(X)=\omega \geq |Y|=l_{u}(Y)$ and $\varphi (Y)=\omega \geq |Y|^{+}=l(Y)$,
so (iv) holds.
\end{proof}

This corollary implies the following results: (1) The Browder Fixed Point
Theorem \ref{Br}, (2) Theorem 7 of Browder \cite{Br}, (3) Theorem 4.3 of
Horvath \cite{Ho1} and Theorem 3.1 of Ding and Tarafdar \cite{DT} for
H-spaces, (4) Browder-type Fixed Point Theorem of Van de Vel \cite[IV.6.28,
p. 506]{VdV}.\ 

We can generalize the definition of convexity by considering a system of
``approximative'' convexity multifunctions converging to $C.$ As a result we
can obtain the results of this paper for such spaces as AES \cite{BD}, ANES 
\cite{GG2}, the ``comb space'', and admissible in the sense of Klee subsets
of topological vector spaces \cite{Had0}.

\section{Appendix: Other Definitions of Generalized Convexity.\label%
{Appendix}}

In this section we provide (without proof) comparison of our definition of
convexity with those due to Horvath, Van de Vel and Michael.

\subsection{H-spaces.\label{Horvath}}

The following notion, originating from the work of Horvath \cite{Ho,Ho1}, is
a generalization of the convex hull in a topological vector space.

\begin{definition}
\label{Hsp}A pair $(Z,\{\Gamma _{A}\})$ will be called an \textit{H-space},
if $Z$ is a topological space and $\{\Gamma _{A}\}$ is a family of
contractible subsets of $Z$ indexed by all finite subsets of $Z$ so that 
\begin{equation*}
\Gamma _{A}\subset \Gamma _{B}\text{ whenever }A\subset B.
\end{equation*}
($(Z,\{\Gamma _{A}\})$ is called a \textit{c-space} \cite{Ho1}). A set $%
A\subset Z$ is called \textit{H-convex} if $\Gamma _{D}\subset A$ for any
finite $D\subset A$, and the \textit{H-convex hull} of a set $A\subset Y$ is
given by 
\begin{equation*}
conv^{\ast }(A)=\cup \{\Gamma _{D}:D\subset A,\ D\text{ is finite}\}.
\end{equation*}
\end{definition}

\begin{proposition}
\label{4.6}Let $(Z,\{\Gamma _{A}\})$ be an H-space$.$ Then there is a global
discrete convexity $(Y,C,Z)$ associated with$\ $some $C:\Delta
(Y)\rightarrow Y,$ $\mathcal{A}=2^{Y}\backslash \{\emptyset \}$, and $conv$
given by: $conv(A)=conv^{\ast }(A),\ A\in \mathcal{A}.$
\end{proposition}

The next definition provides an analogue of local convexity.

\begin{definition}
\label{Hsp1}(1) (Horvath \cite{Ho1}) We say that an H-space $(Y,\{\Gamma
_{A}\})$ is an \textit{l.c.-space} if there is a uniform base $\mathcal{B}$
such as for any $U\in \mathcal{B}$, 
\begin{equation*}
B(A,U)\text{ is an H-convex set whenever }A\subset Y\text{ is H-convex.}
\end{equation*}
A metric space $(Y,d)$ is called a \textit{metric l.c. space} if it is a
c-space and $\forall \varepsilon >0,\ \{y\in Y:d(y,E)<\varepsilon \}$ is an
H-convex set if $E$ is an H-convex set, and open balls are H-convex.

\noindent (2) (Had\v{z}i\'{c} \cite{Ha}) We say that an H-space $(Y,\{\Gamma
_{A}\})$ is of \textit{generalized Zima type }if there is a uniform base $%
\mathcal{B}$ such that for every $U\in \mathcal{B,}$ there exists a $V\in 
\mathcal{B}$ such that for every finite subset $D$ of $Y$ and every H-convex
subset $A$ of $Y$ the following holds: 
\begin{equation*}
A\cap B(z,V)\neq \emptyset \text{ for every }z\in D\Longrightarrow A\cap
B(u,U)\neq \emptyset \text{ for every }u\in \Gamma _{D}.
\end{equation*}
\end{definition}

\begin{proposition}
\label{4.8}Suppose that

(1) $(Y,\{\Gamma _A\})$ is an l.c.-space, or

(2) $(Y,\{\Gamma _{A}\})$ is an H-space of generalized Zima type with
H-convex points.

\noindent Then there is a regular convexity $(Y,C,Y)$ associated with $%
\mathcal{C}=\{$H-convex sets$\}$.
\end{proposition}

\subsection{Van de Vel's Uniform Convex Structures.}

\begin{definition}
\label{VdVdef}\cite[pp. 3 and 304]{VdV} A pair $(Y,\mathcal{C)}$, where $%
\mathcal{C}$ is a family of subsets of $Y$, called \textit{convex sets}, is
called a \textit{uniform convex structure} if

(1) The empty set $\emptyset $ and the universal set $Y$ are in $\mathcal{C}$%
,

(2) $\mathcal{C}$ is stable for intersections, that is, if $\mathcal{D}%
\subset \mathcal{C}$ is nonempty, then $\cap \mathcal{D}$ is in $\mathcal{C}$%
,

(3) $\mathcal{C}$ is stable for nested unions, that is, if $\mathcal{D}%
\subset \mathcal{C}$ is nonempty and totally ordered by inclusion, then $%
\cup \mathcal{D}$ is in $\mathcal{C,}$

(4) there is a uniform base $\mathcal{B}$ such that for each $U\in \mathcal{%
B,}$ there is a $V\in \mathcal{B}$ such that 
\begin{equation*}
\text{for any }A\in \mathcal{C}\text{, }conv(B(A,V))\subset B(A,U),
\end{equation*}
where the \textit{convex hull} $conv$ is defined as: 
\begin{equation*}
conv(A)=\bigcap \{D\in \mathcal{C}:A\subset D\},\ A\subset Y.
\end{equation*}
\end{definition}

\begin{proposition}
\label{VdVConv}If $(Y,\mathcal{C})$ is a uniform convex structure such that
all elements of $\mathcal{C}$ are AR's, then there is a global continuous
convexity $(Y,C,Y)$ associated with $\mathcal{C.}$
\end{proposition}

Van de Vel says that his convex structure satisfies the $S_{4}$\textit{-axiom%
} if two disjoint convex sets can be separated by two convex sets complement
to each other. He calls convex hulls of finite sets \textit{polytopes}. The
statement below follows from Van de Vel's selection theorem. The proposition
below follows from \cite[Theorem 3.17, p. 446]{VdV}.

\begin{proposition}
\label{VdVconv} Let $(Y,\mathcal{C})$ be a metrizable Van de Vel uniform
convex structure satisfying the $S_{4}$-axiom with compact polytopes such
that $\mathcal{C}$ contains only connected sets. Then there is a global
continuous convexity $(Y,C,Y)$ associated with $\mathcal{C.}$
\end{proposition}

\subsection{Michael Convex Structures.\label{Michael}}

\begin{proposition}
\label{MichConv}Let $\{(M_{n},k_{n})\}$ be a Michael's convex structure$.$
Let 
\begin{equation*}
\begin{array}{c}
\mathcal{A}=\{A\subset Y:B(A,W)^{n+1}\subset M_{n}\text{ for all }n\geq 0,\
W\in \mathcal{B}\}, \\ 
conv(A)=\{k_{n}(t,x):x\in A^{n+1},\ t\in \Delta _{n},\ n\geq 0\},\ A\in 
\mathcal{A}
\end{array}
\end{equation*}
(notice that $\mathcal{A}$ can be empty even if $M_{n}$ are not). Then $%
(Y,C,Y)$ is a continuous convexity if $C:Q^{\prime }\rightarrow Y$, where $%
Q^{\prime }=\cup_{A\in \mathcal{A,}W\in \mathcal{B}}\Delta (B(A,W))$
(cf. (c8)), is defined as follows: for any $(d,a)\in Q^{\prime },$%
\begin{equation*}
C(d,a)=k_{n}(t,x),
\end{equation*}
where 
\begin{equation*}
\begin{array}{l}
t=(d_{i_{0}},\ldots ,d_{i_{n}})\text{ and }\ x=(a_{i_{0}},\ldots ,a_{i_{n}}),
\\ 
\{i_{k}:k=0,\ldots ,n\}=\{i\in I:d_{i}\neq 0\},\ i_{0}<\ldots <i_{n},
\end{array}
\end{equation*}
(we assume that the index set $I$ is totally ordered).
\end{proposition}

Let's consider an example of a Michael convex structure: $%
Y=M_{0}=\{0,1\}\subset \mathbf{R},\ M_{n}=\emptyset $ for $n\geq 1,\
k_{0}(t,x)=x$ for all $x\in Y,t\in \lbrack 0,1].$ Let $E=Y=\{0,1\},\
Q_{0}(t,x)=x$ for all $x\in Y,\ t\in \lbrack 0,1]$. Then it is obvious that $%
Q_{1}$ does not exist. More generally, we can equip an $m$-sphere $\mathbf{S}%
^{m},\ m\geq 0$, with a non-trivial Michael convex structure by using its
local Euclidean structure. These examples show that Michael convex
structures are not generalized by convex structures due to Park and Kim \cite
{PK} and Pasicki \cite{Pas}.

\end{document}